\def\zbar{\overline{z}}
\def\tr{{\rm tr}}

\def\IH{{\mathbb H}^3}

\def\IR{{\mathbb R}}

\def\ID{{\mathbb D}}

\def\IS{{\mathbb S}}
\def\IZ{{\mathbb Z}}
\def\IH{{\mathbb H}}
\def\IC{\mathbb C}

\def\oC{\hat{\IC}}

\def\zbar{{\overline{z}}}
\def\hbar{{\overline{h}}}

\documentclass{article}

\usepackage{amsthm, amssymb}
\usepackage{amsmath}
\usepackage{amsfonts}
\usepackage{epsfig,multicol}
\usepackage{pstricks}
\usepackage{graphicx}

\title{Random Ideal Hyperbolic Quadrilaterals, the Cross Ratio Distribution and Punctured Tori.}

\author{Gaven J. Martin  \thanks{Research supported in
part by grants from the New Zealand
Marsden Fund. 
}  }

\date{}

\begin{document}

\maketitle
\newtheorem{theorem}{Theorem}[section]    
                                           
\newtheorem{lemma}[theorem]{Lemma}         
\newtheorem{corollary}[theorem]{Corollary} 
\newtheorem{remark}[theorem]{Remark}       
\newtheorem{definition}[theorem]{Definition}
\newtheorem{conjecture}[theorem]{Conjecture}
\newtheorem{proposition}[theorem]{Proposition}
\newtheorem{example}[theorem]{Example}

\newcommand{\param}{(\gamma,\beta,\beta')}
\newcommand{\parfour}{(\gamma,\beta,-4)}

\numberwithin{equation}{section}

\newcommand{\abs}[1]{\lvert#1\rvert}

\renewcommand{\theequation}{\thetheorem} 
                    
\makeatletter
\let \c@equation=\c@theorem

\begin{abstract}  Earlier work introduced a geometrically natural probability measure on the group of all M\"obius transformations of the hyperbolic plane so as to be able to study ``random'' groups of M\"obius transformations,  and in particular random two-generator groups.  Here we extend these results to consider random punctured tori.  These Riemann surfaces have finite hyperbolic area $2\pi$ and fundamental group the free group of rank 2.  They can be obtained by pairing (identifying) the opposite sides of an ideal hyperbolic quadrilateral.  There is a natural distribution on ideal quadrilateral given by the cross ratio of their vertices.  We identify this distribution and then calculate the distributions of various geometric quantities associated with random punctured tori such as the base of the geodesic length spectrum and the conformal modulus,   along with more subtle things such as the distribution of the distance in Teichm\"uller space to the central ``square" punctured torus.   \end{abstract}

\maketitle

\section{Introduction.}  Earlier work \cite{MO}  introduced a geometrically natural probability measure on the group of all M\"obius transformations,  or isometries,  of the hyperbolic plane with an aim  to study ``random'' groups of M\"obius transformations,  and in particular random two-generator groups - groups where the generators are selected randomly from this distribution.  We determined the basic statistics of these groups and discussed the probability that a random group was discrete.  We gave further explicit results in the case of groups of isometries of hyperbolic $3$-space generated by two randomly selected parabolic elements,  \cite{MOY}.  When such a  group is discrete the boundary of the orbit space is generically a single triply punctured sphere, and thus having no conformal moduli.

\medskip

Here we extend these results to consider a case where conformal moduli and Teichm\"uller spaces appear - albeit in the simplest possible case of random punctured tori.  These Riemann surfaces have area $2\pi$ when endowed with the complete hyperbolic metric of  curvature $-1$,  and fundamental group isomorphic to the free group on two generators.  They are naturally obtained by pairing opposite sides of an ideal hyperbolic quadrilateral $Q$ in the hyperbolic plane.   There is a natural distribution on ideal quadrilaterals given by the cross ratio of their vertices - denoted $[Q]$,  (precise definitions below).  Thus our first task is to identify the distribution of cross ratios for four points randomly (and uniformly) selected on the circle.  This is Theorem \ref{crdist}.  

Next,  given a random ideal quadrilateral $Q$ we can topologically identify (pair) its sides. It is not difficult to see that if we do this randomly $\frac{2}{3}$ of the time we get a hyperbolic triply punctured sphere $S^{2}_{3}$,  and $\frac{1}{3}$ of the time we get a once punctured torus $T^{2}_{*}$.  The dichotomy being determined by whether we pair adjacent or opposite sides. The surface $S^{2}_{3}$ has no conformal invariants at all since all triply punctured spheres are isometric.  We therefore explore the latter case.  We seek the distributions and expected values of geometric invariants of $T^{2}_{*}$ such as the length of the shortest geodesic $\ell_{[Q]}$,  where we obtain answers in closed form,  and the conformal modulus of $T^{2}_{*}$ denoted $m(Q)$.  This latter association $[Q]\leftrightarrow m(Q)$ is a difficult classical problem.  There is no closed form solution, and  the answer lies beyond the theory of special functions.  Thus we present numerical investigations backed up with rigorous  asymptotic calculations.  For instance we give a computational description of the probability distribution function for the distance in the Teichm\"uller metric of a random punctured torus to the origin - the ``square'' punctured torus.

\bigskip

\noindent {\bf Remark.} We note that in the particular case of $T^{2}_{*}$ there is another approach, \cite{MarT}.  Since the moduli space is know to us as a finite area Riemann surface,  we can put the uniform distribution on it, and from this calculate the distribution of other quantities such as distance to the square punctured tori.  However this approach is not especially useful as almost no other distributions of  invariants  can be readily calculated,  for instance even the length of the shortest geodesic can't be given in a closed form such as Theorem \ref{5.7} later.  Also the p.d.f. for the Teichm\"uller distance obtained from this approach has a couple of quite strong singularities. Finally,  the approach  found in this paper offers the opportunity to consider every finite area non-compact surface as these can be obtained by side pairings of ideal polyhedra,  and we may select the vertices of these polyhedra randomly and uniformly in $\IS$,  whereas the higher dimensional Teichm\"uller spaces are far more complicated objects to deal with metrically.

\bigskip

\noindent {\bf Thanks.}  The author wishes to thank Alex Eremenko for guiding him through his article \cite{Eremenko} and proving some assistance with the code for the computational results of \S 7 \& 8.

\section{Ideal Hyperbolic Quadrilaterals.}

If we identify the hyperbolic plane $\IH^2$ with the unit disk $\ID$ and hyperbolic metric $ds=\frac{|dz|}{1-|z|^2}$, then an ideal hyperbolic quadrilateral is a quadrilateral all of whose sides are hyperbolic lines and all the vertices lie on the circle at infinity - identified as the boundary $\IS=\partial\ID$.   Equivalently this is the region bounded by the hyperbolic convex hull of four points in $\IS$.  There is a single conformal invariant of an ideal quadrilateral - the cross ratio of the vertices $z_1,z_2,z_3,z_4$ (see (\ref{crdef}) below) up to the permutations of the vertices as described at (\ref{S4}). Other conformal invariants - such as the hyperbolic distance between the sides - are therefore functions of this cross ratio. One can easily see the geometric action of  $S_4$ on the symmetric quadrilateral with vertices $\{\pm 1, \pm i \}$.
  \begin{center}
\scalebox{0.55}{\includegraphics*[viewport=20 560 620 830]{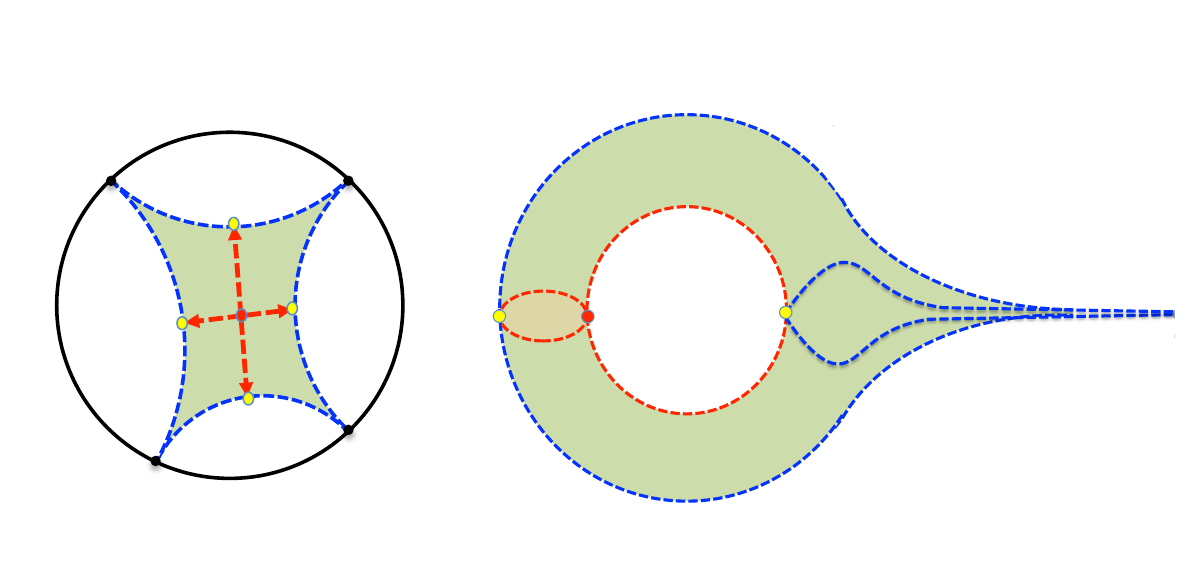}} \\
\end{center}
 {\bf Figure 1.} {\em  Left: An ideal quadrilateral in the hyperbolic disk $\ID$.  Right:  A rectangular punctured torus obtained by identifying opposite edges.}

\medskip

Conversely,  as noted, the hyperbolically convex set generated by four distinct points on the circle defines a unique quadrilateral.  It is an easy exercise to see that two such quadrilaterals  $Q_z={\rm Convex Hull}\big(\{z_1,z_2,z_3,z_4\}\big)$ and $Q_w={\rm Convex Hull}\big(\{w_1,w_2,w_3,w_4\}\big)$ are isometric if and only if the cross ratios $\lambda_z=[z_1,z_2,z_3,z_4]$ and $\lambda_w=[w_1,w_2,w_3,w_4]$ satisfy
\[ 
  \{ \lambda_z,   1-\lambda_z, \frac{ \lambda_z}{\lambda_z-1}, \frac{1}{\lambda_z}, \frac{1}{1-\lambda_z}, \frac{\lambda_z\-1}{\lambda_z}\} \cap     \{ \lambda_w,   1-\lambda_w, \frac{ \lambda_w}{\lambda_w-1}, \frac{1}{\lambda_w}, \frac{1}{1-\lambda_w}, \frac{\lambda_w\-1}{\lambda_w}\}  \neq \emptyset.
\]

\section{The Cross Ratio}

Given four points $z_1,z_2,z_3,z_4 \in \oC$ we defined their {\em cross ratio} as
\begin{equation}\label{crdef}
[z_1,z_2,z_3,z_4] = \frac{(z_1-z_3)(z_2-z_4)}{(z_1-z_2)(z_3-z_4)}.
\end{equation} 
 The permutation group $S_4$ on the set $\{1,2,3,4\}$ acts on this cross ratio as
 \[ [z_1,z_2,z_3,z_4] \to [z_1,z_2,z_3,z_4]_\sigma = [z_{\sigma(1)},z_{\sigma(2)},z_{\sigma(3)},z_{\sigma(4)}],  \hskip10pt \sigma\in S_4. \]
In this way the group $S_4$ effects the transformations of $\lambda=[z_1,z_2,z_3,z_4] $ by  
\begin{equation}\label{S4}
 \lambda\mapsto\lambda,  \;\; 1-\lambda, \;\;   \lambda/(\lambda-1), \;\;  1/\lambda, \;\;  1/(1-\lambda), \;\;  1-1/\lambda.
\end{equation}
This is all carefully explained in \cite[Theorem 4.4.1]{Beardon}. Therefore given four points randomly selected from the uniform distribution on the circle the probability distribution function (p.d.f.) must reflect these symmetries.  We may select four such points as
 \[ z_j= e^{i \theta_j}, \hskip10pt j=1,2,3,4,   \hskip10pt \theta_j\in_u [0,2\pi],  \]
 the last term here meaning we select $\theta_j$ uniformly in $[0,2\pi]$. 
 
\medskip

We then calculate 
\begin{eqnarray*} 
 \frac{(e^{i \theta_1}-e^{i \theta_3})(e^{i \theta_2}-e^{i \theta_4})}{(e^{i \theta_1}-e^{i \theta_2})(e^{i \theta_3}-e^{i \theta_4})} & = &  
  \frac{ \sin \left(\frac{\theta_1 -\theta_3 }{2}\right) \sin \left(\frac{\theta_2 -\theta_4 }{2}\right)}{ \sin \left(\frac{\theta_1 -\theta_2 }{2}\right) \sin \left(\frac{\theta_3 -\theta_4 }{2}\right)},
 \end{eqnarray*}
 and so we see that we may as well select the angles uniformly  $\theta_j\in_u[0,\pi]$ and consider the distribution of
 \[ \frac{ \sin \left(\theta_1 -\theta_3\right) \sin \left(\theta_2 -\theta_4 \right)}{ \sin \left(\theta_1 -\theta_2 \right) \sin \left(\theta_3 -\theta_4 \right)} . \]

Further since rotation leaves both Lebesgue measure and cross ratio invariant we may as well assume $\theta_4=0$. More carefully,  observe that $(\theta_1 -\theta_4)$ mod $\pi$,  $(\theta_2 -\theta_4)$ mod $\pi$ and $(\theta_3 -\theta_4)$ mod $\pi$ are uniformly distributed in $[0,\pi]$ so we can use these as variables.  Now
 \begin{eqnarray*}
\lambda  & = &   \frac{ \sin \left(\theta_1 -\theta_3\right) \sin \left(\theta_2\right)}{ \sin \left(\theta_1 -\theta_2 \right) \sin \left(\theta_3\right)}  = \frac{ \big( \sin(\theta_1)\cos(\theta_3)-\sin(\theta_3)\cos(\theta_1) \Big) \sin(\theta_2)}{ \big(\sin(\theta_1)\cos(\theta_2)-\sin(\theta_2)\cos(\theta_1)\big)   \sin \left(\theta_3\right)} \\
 & = & \frac{  \cot(\theta_3)-\cot(\theta_1)   }{  \cot(\theta_2)-\cot(\theta_1)    } .
 \end{eqnarray*}
 
 Next,  the mapping 
 \begin{equation}
 \IR^3 \ni (x,y,z)\mapsto (\cot^{-1}(x),\cot^{-1}(y),\cot^{-1}(z)) \in [0,\pi]^3, 
 \end{equation}
 is an orientation reversing diffeomorphism.  It's Jacobian determinant is 
\begin{equation} \label{jac}{\bf J}(x,y,z)  =  - \frac{1}{1+x^2} \; \frac{1}{1+y^2} \; \frac{1}{1+z^2} . \end{equation}
Therefore the change of variables formula gives us a simply expression for the cumulative distribution of the cross ratio.
\begin{eqnarray}\label{cdf}
\nonumber F(r)={\text {Pr}}\{\lambda < r \} & = & \frac{1}{\pi^3} \; \int_{ \{(\theta_1,\theta_2,\theta_3)\in [0,\pi]^3 :  \frac{  \cot(\theta_3)-\cot(\theta_1)   }{  \cot(\theta_2)-\cot(\theta_1)    }  < \lambda \} } d\theta_1d\theta_2d\theta_3 \\
& = &  \frac{1}{\pi^3} \; \int_{ \{(x,y,z)\in\IR^3 :  \frac{ x-y   }{ x-z  }  < \lambda \} }  {\bf J}(x,y,z)  \; dxdydz.
\end{eqnarray}
 
 \section{The Cross Ratio Distribution.}
 
 We wish to calculate the derivative with respect to $r$ of the cumulative distribution function $F(r)$ at (\ref{cdf}) where the Jacobian ${\bf J}$ is defined at (\ref{jac}).

\medskip

 \subsection{$r$ negative.}    We begin this calculation under the assumption $r<0$.  Then
  \begin{eqnarray*}
\pi^3\; F(r) & =  & \int_{\{(x,y,z): \frac{x-y}{x-z} < r \}} {\bf J}(x,y,z) \; dxdydz \\
& =  & \int_{\{  x-y  < r(x-z) \} \cap \{x>z\}} {\bf J}  \; + \int_{\{ x-y > r(x-z) \} \cap \{x<z\}} {\bf J}  \\
& =  & \int_{\{  z<x  < \frac{y-rz}{1-r} \}} {\bf J}  \; + \int_{\{  \frac{y-rz}{1-r} < x< z \} } {\bf J}  \\
& =  & \int_{-\infty}^{\infty} \; \int_{z}^{\infty} \;  \int_{z}^{\frac{y-rz}{1-r} } {\bf J}  \; dx \; dy \; dz + \int_{-\infty}^{\infty} \; \int^{z}_{-\infty} \;  \int_{\frac{y-rz}{1-r}}^{z}  {\bf J} \; dx \; dy \; dz  .
\end{eqnarray*}
We can now differentiate this with respect to $r$ to see that
 \begin{eqnarray*}
\pi^3\;F'(r) & =  &  
 \int_{-\infty}^{\infty} \; \int_{z}^{\infty} \; \frac{y-z}{(1-r)^2} \;  {\bf J}\Big|_{x= \frac{y-rz}{1-r}}   \; dy \; dz - \int_{-\infty}^{\infty} \; \int^{z}_{-\infty}  \frac{y-z}{(1-r)^2} \;  {\bf J}\Big|_{x= \frac{y-rz}{1-r}}    \; dy \; dz  \\
 & =  &  
 \int_{-\infty}^{\infty} \; \int_{-\infty}^{\infty} \; \frac{|y-z|}{((1-r)^2 +( y-rz)^2)(1+y^2)(1+z^2) }   \; dy dz. \end{eqnarray*}
 Now
  \begin{eqnarray*}
  \int^{z}_{-\infty} \; \frac{z-y}{((1-r)^2 +( y-rz)^2)(1+y^2)(1+z^2) }  \; dy 
& = & \frac{r z \left(2 \tan ^{-1}(z)+\pi \right)+(r-2) \log (1-r)}{r \left(z^2+1\right) \left(r \left(r z^2+r-4\right)+4\right)}.
 \end{eqnarray*}
 We therefore want to evaluate
 \[  {\bf I}_1 =\int_{-\infty}^{\infty} \; \frac{r z \left(2 \tan ^{-1}(z)+\pi \right)+(r-2) \log (1-r)}{r \left(z^2+1\right) \left(r \left(r z^2+r-4\right)+4\right)} \; dz .\]
 We break this up and consider three terms here.  First, 
   \begin{eqnarray*}
\int^{\infty}_{-\infty} \; \frac{r z  \pi }{r \left(z^2+1\right) \left(r \left(r z^2+r-4\right)+4\right)} = 0.
 \end{eqnarray*}
 Then 
   \begin{eqnarray*}
\int^{\infty}_{-\infty} \;  \frac{(r-2) \log (1-r)}{r \left(z^2+1\right) \left(r \left(r z^2+r-4\right)+4\right)} \; dz  = \frac{\pi  \log (1-r)}{2 (r-1) r}.
 \end{eqnarray*} 
 Finally the more challenging integral.
    \begin{eqnarray*}
\int_{-\infty}^{\infty} \; \frac{2 z  \tan ^{-1}(z)}{ \left(z^2+1\right) \left(r \left(r z^2+r-4\right)+4\right)} \; dz & = & \int^{\pi/2}_{-\pi/2} \; \frac{ \tan(w)   w \; dw}{ r^2\sec^2(w)+4-4r}  \\
& = & \frac{\pi  \log \left(1-\frac{1}{r}\right)}{2 (1-r)}.
 \end{eqnarray*} 
 We therefore obtain
 \begin{eqnarray*}
 {\bf I}_1 & = & \frac{\pi  \log (1-r)}{2 (r-1) r} + \frac{\pi  \log \left(1-\frac{1}{r}\right)}{2 (1-r)} .
 \end{eqnarray*}
Similary,  we must calculate 
   \begin{eqnarray*}
\lefteqn{   \int^{\infty}_{z} \; \frac{y-z}{((1-r)^2 +( y-rz)^2)(1+y^2)(1+z^2) }  \; dy } \\
& = & -\frac{-(1-2 r) \log \left(\frac{r-1}{r}\right)+\frac{\pi  z}{2}+z \tan ^{-1}\left(\frac{1}{z}\right)-z \tan ^{-1}(z)}{\left(z^2+1\right) \left(4 r^2-4 r+z^2+1\right)}.
 \end{eqnarray*}
 When we integrate from $-\infty$ to $\infty$ the term involving $\frac{\pi  z}{2}$ vanishes and we are left with ${\bf I}_1$ again.  We have proved the following theorem once we include the symmetries given by the $S_4$ action described above to extend the range of the p.d.f. to positive values of $r$ as well.
 
 \begin{theorem}\label{crdist}  Let $\lambda=[z_1,z_2,z_3,z_4]$ be the random variable where $z_i$ are selected uniformly in the circle,  $z_i \in_{u} \IS^1$.  Then the probability distribution function   $\lambda(r)$ is
 \begin{eqnarray}
\frac{1}{\pi^2}\left(  \frac{  \log (1-r)}{ (r-1) r} + \frac{  \log \left(1-\frac{1}{r}\right)}{(1-r)}  \right), && \mbox{if $r\leq 0$}. \\
\frac{1}{\pi^2}\left(\frac{ \log \big(\frac{1}{1-r}\big)}{r}+\frac{  \log \left(\frac{1}{r}\right)}{(1-r)}\right) , && \mbox{if $0\leq r\leq 1$}. \\
\frac{1}{\pi^2}\left( \frac{  \log (r)}{(r-1) r} + \frac{ \log \left( \frac{r}{r-1}\right)}{r}  \right) ,&& \mbox{if $r\geq 1$}.
 \end{eqnarray}
 \end{theorem}
   \begin{center}
\scalebox{0.6}{\includegraphics*[viewport=75 380 760 700]{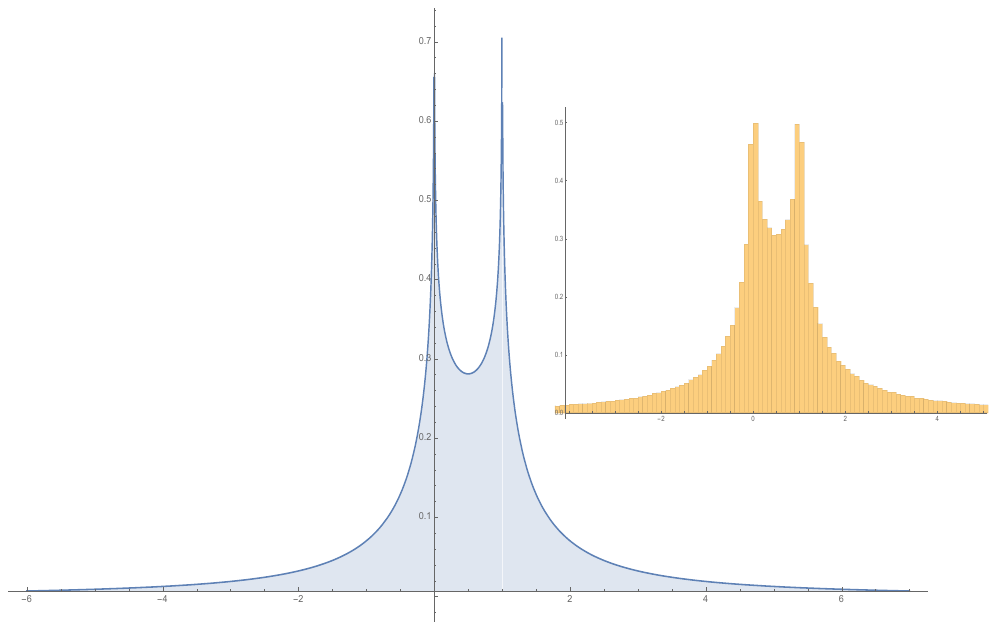}} \\
{\em The probablity distribution function for the cross ratio of four random points on the circle.  Inset a p.d.f histogram from $10^7$ random quadruple cross ratios.}
\end{center}

\subsection{Random ideal quadrilaterals.}
We have defined a random ideal quadrilateral to be the hyperbolic convex hull of four randomly (and uniformly) selected points on the circle.  We can unambiguously assign a cross ratio to a quadrilateral as follows.  We cyclically order the vertices counterclockwise.  The quadrilateral is M\"obius equivalent to that with vertices $\{e^{i\theta},-e^{-i\theta},-e^{i\theta},e^{-i\theta}\}$ where we have the additional freedom to choose $\theta\in[0,\pi/4]$ yielding cross ratio
\[ [e^{i\theta},-e^{-i\theta},-e^{i\theta},e^{-i\theta}]=\frac{1}{\cos^2 \theta} \geq 2. \]
For a random quadrilateral we define $[Q]$ as the unique cross ratio associated with the vertices and having value at least $2$.  The p.d.f. for $[Q]\geq 2 $ is then
\begin{equation}\label{Qdef} X_{[Q]} = \frac{6}{\pi^2}\left[ \frac{  \log (r)}{(r-1) r} + \frac{ \log \left( \frac{r}{r-1}\right)}{r}  \right] , \hskip5pt r\geq 2. \end{equation}
The closed form for the cumulative distribution is
\[  \frac{6}{\pi^2} \left[\text{Li}_2\Big(\frac{r}{r-1}\Big)-\text{Li}_2(r) \right]+\frac{3}{\pi ^2}   \log (r-1) \; \log \Big(\frac{r-1}{r^2}\Big).  \]

 \subsection{Median}
We find the median value of $X_{[Q]}$ to be the solution $\gamma$  to the equation
\[\text{Li}_2\left(\frac{\gamma}{\gamma-1}\right)-\text{Li}_2(\gamma)+\frac{1}{2} \log (\gamma-1) (\log (\gamma-1)-2 \log (\gamma)-2 i \pi )=-\frac{\pi ^2}{12} ,\]
and hence $\gamma \approx 4.6883\ldots $.

There are no higher moments for $X_{[Q]}$ as
\begin{equation}
\frac{\pi^2}{6} X_{[Q]}(r) =\frac{ \log (r)+1}{ r^2}+\frac{ \log (r)+\frac{1}{2}}{r^3}+O(r^{-4}),  \hskip10pt r >> 1.
\end{equation}

 \subsection{Normalised cross ratio.}  If we had normalised a cross ratio via M\"obius transformations,  then we can chose $z_1=1$, $z_2=i$, $z_3=-1$ with free variable $z_4=e^{i\theta}$ with $\theta\in_u [0,2\pi]$.  The formula for the cross ratio is then
 \[ \lambda^*= [i,1,-1,z_4] =   \frac{(i+1)(1-e^{i\theta})}{(1-i)(1+e^{i\theta})}  = \tan \theta/2 \]
and the cumulative distribution is found from the fact that
 \[ \big|\{\theta\in [0,\pi]: \tan(\theta) < r \} \big|= \left\{ \begin{array}{cc} \tan^{-1}(r)-\frac{\pi}{2}  & \mbox{if $r<0$ } \\ \tan^{-1}(r)+\frac{\pi}{2}  & \mbox{if $r>0$ }\end{array} \right. \]
 Differentiating this gives us the p.d.f. for $\lambda^*$ as $X^*(r)=\frac{1}{\pi} \; \frac{1}{1+r^2}$.
 
 \section{Cross ratios and moduli of punctured tori.}

A punctured torus $T^{2}_{*}$  admits a covering of the form  
\begin{equation}
T^{2}_{*} \approx (\IC \setminus \Lambda(0))/\Lambda,
\end{equation} 
where  $\Lambda$ is a group of translations of $\IC$,
\begin{equation}
\Lambda = \Lambda_{\alpha,\beta} =  \{z\mapsto z+ m\alpha + n \beta: m,n\in \IZ\},
\end{equation}
and $\alpha,\beta\in \IC\setminus\{0\}$,  $\alpha/\beta\not\in\IR$,  are periods.  When $\IC\setminus \Lambda(0)$ is equipped with the complete hyperbolic metric, this induces a hyperbolic metric on $T^{2}_{*} $ and in fact every conformal equivalence class of punctured tori can be obtained this way. Further,  the $SL(2,\IZ)$ action on the set of lattices allows us to normalise $\Lambda$.  We can assume that $\alpha=1$ and that $\beta$ lies in the fundamental domain $G$ for the modular group acting on the upper half-plane,  $G=\{z\in\IC: \Im m(z)>0, |z|\geq 1, -\frac{1}{2}< \Re e(z) \leq \frac{1}{2} \}$.  The lattice $\Lambda$ will be rectangular if $\beta$ is purely imaginary.  However,  to facilitate later computations we normalise a rectangular lattice as follows.   We fix $\beta=i$ and $\Lambda_\alpha=\Lambda_{\alpha,i}$,  $\alpha \geq 1$ and in this case we say the punctured torus $T^{2}_{*}(\alpha)=  (\IC \setminus \Lambda_{\alpha}(0))/\Lambda_{\alpha}$ is rectangular.  Rectangular punctured tori arise as described in Figure 1 from the identification of opposite sides of an ideal hyperbolic quadrilateral when we ensure the closest points of each side pair are identified.

\subsection{Punctured torus groups.}

The fundamental group of a punctured torus is free on two generators.  However there is more structure when we consider a representation as a Fuchsian group.  There we have the geometric presentation $\langle a,b:[a,b]^\infty=1\rangle$,  where the presentation means that the multiplicative commutator is parabolic,  and that $a$ and $b$ are hyperbolic.  This group can be represented in $PSL(2,\IC)$ by
\[ A=\pm \left(
\begin{array}{cc}
 \sqrt{\frac{1}{r^2}+1} & \frac{1}{r} \\
 \frac{1}{r} & \sqrt{\frac{1}{r^2}+1} \\
\end{array}
\right),  \hskip10pt B= \pm \left(
\begin{array}{cc}
 \sqrt{\frac{1}{s^2}+1} & \frac{i}{s} \\
 -\frac{i}{s} & \sqrt{\frac{1}{s^2}+1} \\
\end{array}
\right), \]
and the commutator is parabolic if and only if $rs=1$,  and then simplifies to 
\[ [A,B] = \left(
\begin{array}{cc}
 -2 i \sqrt{r^2+1} \sqrt{s^2+1}-1 & 2 i \sqrt{s^2+1}-2 \sqrt{r^2+1} \\
 -2 \sqrt{r^2+1}-2 i \sqrt{s^2+1} & 2 i \sqrt{r^2+1} \sqrt{s^2+1}-1 \\
\end{array}
\right).\]
The matrices $A$ and $B$ represent the M\"obius transformations
\begin{equation}\label{Af} f(z)=\frac{ 
 \sqrt{ {r^2}+1} \; z + 1}
 { z +\sqrt{ r^2+1}}, \hskip10pt   g(z)=\frac{  \sqrt{ {s^2}+1} z + i}{
 -i z + \sqrt{ {s^2}+1}}.
 \end{equation}
   The fixed points of $f$ are $\pm1$ and those of $g$ are $\pm i$.  Both $f$ and $g$ setwise fix the unit circle and act as isometrices of the hyperbolic disk $\ID$. The mappings $f$ and $f^{-1}$ pair their isometric circles ${\cal C}(f) = \{z:|f'(z)|=1\}$,  ${\cal C}(f^{-1}) = \{z:|(f^{-1})'(z)|=1\}$.
 \[{\cal C}(f) = \{z:|z + \sqrt{ {r^2}+1}|=r \}, \hskip10pt {\cal C}(f^{-1}) = \{z:|z - \sqrt{ {r^2}+1}|=r \},\]
  \[{\cal C}(g) = \{z:|z + i\sqrt{ {s^2}+1}|=s \}, \hskip10pt {\cal C}(g^{-1}) = \{z:|z - i\sqrt{ {s^2}+1}|=s \}.\]
 When $rs=1$ the circles ${\cal C}(f)$ and ${\cal C}(g)$ are tangent,  and similarly ${\cal C}(f^{-1})$ and ${\cal C}(g^{-1})$ are tangent. Then $\ID \cap \{{\cal C}(f^{\pm 1}) \cup {\cal C}(g^{\pm 1})\}$ bounds an ideal hyperbolic quadrilateral $Q_{r,s}$ in $\ID$.  It follows that $Q_{r,s}$ is a fundamental polygon for the group $\langle f,g \rangle$ acting on $\ID$.  A straightforward geometric calculation reveals the following lemma.
 \begin{lemma}  The cross ratio (up to $S_3$ symmetry) of the vertices of $Q=Q_{r,1/r}$ is
 \begin{equation}
 [Q]= 1+r^2.
 \end{equation}
 and $r,\frac{1}{r}$ are the radii of the isometric circles of generators.
 \end{lemma}
 Next,  we can use the cross ratio to determine the hyperbolic length of the common perpendiculars of the edges of $Q$.  For $t>1$, the four points $-t,-1,1,t \in \partial \IH$ have cross ratio 
 \[ [-t,-1,1,t]=\Big(\frac{1+t}{1-t}\Big)^2. \]
 The hyperbolic distance between the hyperbolic lines with endpoints $-t,t$ and endpoints $-1,1$ is $\log t$.  This shows us that the largest length of the two common perpendiculars between opposite sides has length $\ell_r$ where
 \[ \ell_r = \log \frac{\sqrt{[Q]}+1}{ \sqrt{[Q]}-1},  \hskip10pt [Q] = \coth^2 \Big(\frac{\ell_r}{2}\Big) , \] 
 so $1+r^2 = \coth^2 \Big(\frac{\ell}{2}\Big)$ and  $ r = \sinh \Big(\frac{\ell_r}{2}\Big)$.  Symmetry shows us that $ s = \sinh \Big(\frac{\ell_s}{2}\Big)$ and $1=\sinh(\ell_r)\sinh(\ell_s)$,  all formulas which can be found in \cite{Beardon}.   Since we have the probability distribution for $[Q]$ reported at (\ref{Qdef}) we can compute the distributions for all these quantities.
 
 \begin{theorem} \label{5.7} Let the rectangular punctured torus $T^{2}_{*}$ arise from the side pairings of a random ideal quadrilateral $Q$.  Then  the shortest geodesic has p.d.f.
 \[ X_\ell = \frac{6}{\pi ^2} {\rm csch}(\ell) \Big[4 \log \cosh \frac{\ell}{2}  +2(\cosh (\ell)-1) \log \coth \frac{\ell}{2}  \Big], \hskip5pt 0< \ell \leq \log \frac{\sqrt{2}+1}{\sqrt{2}-1}  \]
The expected value of the shortest geodesic is 
\begin{eqnarray*} E[\ell] & = &  \frac{6}{\pi^2} \int_{0}^{\log(3+\sqrt{2})} {\rm csch}(\ell) \Big[4 \log  \cosh   \frac{\ell}{2}  +2(\cosh (\ell)-1) \log  \coth ^2 \frac{\ell}{2} \Big] \; \ell \; d\ell \\ & \approx & 0.984154 \ldots \end{eqnarray*}
and the median value is $0.99929\ldots$.
 \end{theorem}
\begin{center}
\scalebox{0.75}{\includegraphics*[viewport=60 500 760 760]{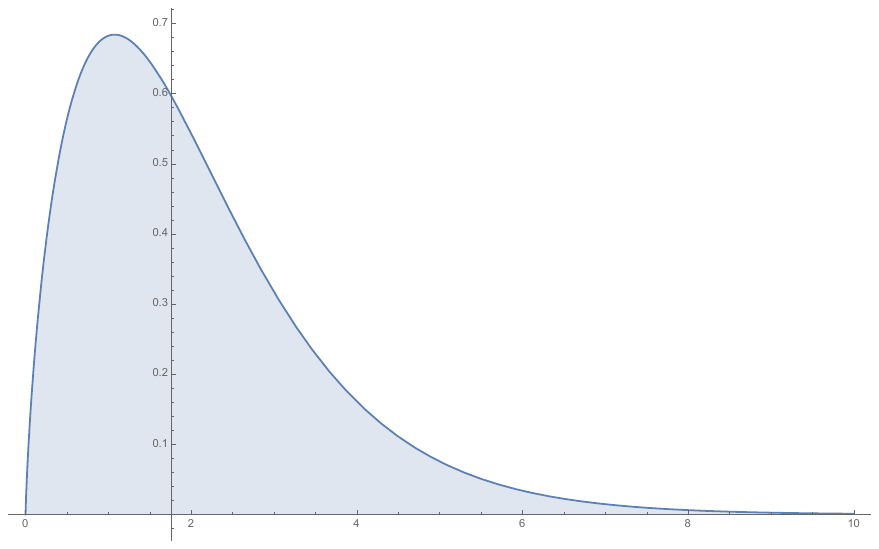}} \\
{\em Length p.d.f. : Shortest geodesic to the left of the vertical axis (at $\log \frac{\sqrt{2}+1}{\sqrt{2}-1}$),  and the length of its dual to the right.}
\end{center}
 \medskip
 
 \subsection{Nonrectangular tori.}  Rectangular tori have opposite sides of a fundamental quadrilateral identified so that the closest points of opposite sides are also identified,  this forces the commutator of the side pairings to be parabolic.  If $Q$ is a quadrilateral with side pairings from the isometric circles of  $f^{\pm1}$ and $g^{\pm1}$,   and   if $h_1$ is hyperbolic and has axis ${\cal C}(f)\cap\ID$ and $h_2$ is hyperbolic and has axis ${\cal C}(g)\cap\ID$,  then $(f\circ h_1)^{\pm 1}$ and $(g\circ h_2)^{\pm 1}$ pair opposite edges of the same quadrilateral (the axes of $h_i$,  $i=1,2$ are adjacent edges of the quadrilateral).  The Poincar\'e Polygon Theorem,  or simple observation shows that $\langle  f\circ h_1, g\circ h_2 \rangle$ has quotient a once punctured torus as the group action wraps up the quadrilateral.   Further all possible conformal classes of once punctured tori can be constructed this way (one needs to observe that every punctured torus has an ideal decomposition by geodesics as a quadrilateral with a single vertex at the cusp,  as per Figure 1.)  
 
Now $f$ has as its axis the common perpendicular between the isometric circles of $f$,  so if we let $\sqrt{f}$ be the hyperbolic with the same axis as $f$ but half the translation length,  then $h_1=\sqrt{f}^{-1}\circ \tilde{g} \circ \sqrt{f}$,  where $g$ and $\tilde{g}$ have the same axis and therefore commute.  Also,
 \[ f \circ h_1 = f\circ \sqrt{f}^{-1}\circ \tilde{g} \circ \sqrt{f} = \sqrt{f} \circ \tilde{g} \circ \sqrt{f},  \]
 and similarly $g\circ h_2 = \sqrt{g} \circ \tilde{f} \circ \sqrt{g} $ with $f$ and $\tilde{f}$ sharing an axis.  We have
 \[ \tilde{f} \sim \tilde{A} =\pm \left(
\begin{array}{cc}
 \sqrt{\frac{1}{\lambda^2}+1} & \frac{1}{\lambda} \\
 \frac{1}{\lambda} & \sqrt{\frac{1}{\lambda^2}+1} \\
\end{array}
\right),  \hskip10pt \tilde{g} \sim \tilde{B} = \pm \left(
\begin{array}{cc}
 \sqrt{\frac{1}{\mu^2}+1} & \frac{i}{\mu} \\
 -\frac{i}{\mu} & \sqrt{\frac{1}{\mu^2}+1} \\
\end{array}
\right),\]
and 
\[ \sqrt{f} \sim \frac{1}{\sqrt{2}}\; \left(
\begin{array}{cc}
1/{\sqrt{r \sqrt{r^2+1}-r^2}} & \sqrt{\sqrt{\frac{1}{r^2}+1}-1} \\
 \sqrt{\sqrt{\frac{1}{r^2}+1}-1} & 1/{\sqrt{r \sqrt{r^2+1}-r^2}} \\
\end{array}
\right), \]
\[ \sqrt{g}  \sim \frac{1}{\sqrt{2}} \left(
\begin{array}{cc}
1/{\sqrt{s \left(\sqrt{s^2+1}-s\right)}} & i \sqrt{\sqrt{\frac{1}{s^2}+1}-1} \\
 -i \sqrt{\sqrt{\frac{1}{s^2}+1}-1} &1/{\sqrt{s \left(\sqrt{s^2+1}-s\right)}} \\
\end{array}
\right). \]
 In order that the isometric circles of $g$ and $g$ form a quadrilateral we have
 $rs = 1$.
It is a lengthy calculation to see that with $u=\sqrt{f} \circ \tilde{g} \circ \sqrt{f}$, $v= \sqrt{g} \circ \tilde{f} \circ \sqrt{g} $ and with $rs=1$ we have that $u$ and $v$ pair opposite sides of the ideal quadrilateral formed from the isometric circles of $f$ and $g$ and that also
\begin{equation}
\tr [u, v] -2 = -\frac{4 \left(\lambda ^2 \left(\mu ^2+1\right)-2 \sqrt{\left(\lambda ^2+1\right) \left(\mu ^2+1\right)}+\mu ^2+2\right)}{\lambda ^2 \mu ^2}.
\end{equation}
 Then $[u,v]$ is parabolic if and only if $\lambda=\mu>0$.  We replace $\lambda$ by $1/\lambda$ to simplify notation and write
 \[ u = \left(
\begin{array}{cc}
 \sqrt{r^2+1}\sqrt{1+ {\lambda ^2}} &  {\lambda }+i r \sqrt{ 1+\lambda ^2} \\
 {\lambda }-i r \sqrt{ 1+\lambda ^2} & \sqrt{r^2+1}\sqrt{1+ {\lambda ^2}}  \\
\end{array}
\right), \]
\[ v= \left(
\begin{array}{cc}
 \sqrt{ 1+\frac{1}{r^2}}\sqrt{1+ {\lambda ^2} } & \frac{1}{r} \sqrt{ \lambda ^2 +1}+ {i}{\lambda } \\
\frac{1}{r} \sqrt{ \lambda ^2+1}- {i}{\lambda } & \sqrt{ 1+\frac{1}{r^2}}\sqrt{1+ {\lambda ^2} } \\
\end{array}
\right).\]
Then we obtain
\begin{eqnarray*}
\tr^2(u)-4 & = &  4(r^2+1)(\lambda ^2+1) -4 = 4 \sinh^2 \frac{\tau_u}{2},\\
 \tr^2(v)-4 & = &  4\Big(\frac{1}{r^2}+1\Big)(\lambda ^2+1) -4 = 4 \sinh^2 \frac{\tau_v}{2}.
\end{eqnarray*}
The fixed points of $u$ and of $v$ are antipodal,
\[ {\rm fix}(u)=\pm \sqrt{\frac{\lambda +ir \sqrt{\lambda ^2+1} }{\lambda -ir \sqrt{\lambda ^2+1} }}, \;\;\;\; {\rm fix}(u)=\pm  \sqrt{ \frac{ \sqrt{\lambda ^2+1}+i \lambda  r}{\sqrt{\lambda ^2+1}-i \lambda  r} }.\]
The axes of $u$ and $v$ meet at the origin at angle $\theta$ where $\sin(\theta) \sinh \frac{\tau_u}{2} \sinh \frac{\tau_v}{2} = 1$.
 \begin{theorem}  Let $\Sigma=T^{2}_{*}$ be a once punctured torus.  Then there is a geodesic ideal quadrilateral $Q \subset \Sigma$,  $\bar Q= \Sigma$ such that the dual geodesics to the sides of $\Sigma$ have length $\ell_1$ and $\ell_2$ and meet once at angle $\theta\in[0,\pi/2]$ with \[ \sin(\theta) \sinh \frac{\ell_1}{2} \sinh \frac{\ell_2}{2} = 1. \]
 The ideal quadrilateral $Q$ is conformally equivalent to an ideal quadrilateral in $\ID$ whose vertex cross ratio is
 \[ [Q]=1+\frac{\cosh^2 \frac{\ell_1}{2}}{\cosh^2 \frac{\ell_2}{2}} . \]
 \end{theorem}
 In light of what we have found above a natural way to put a distribution on  random punctured tori (not just the rectangular ones)  is to choose two positive real numbers $x_\Sigma$ and $y_\Sigma$ from the distribution of Theorem \ref{5.7}.  That is
  \[ {\bf X}[x] = \frac{3}{\pi ^2} {\rm csch}(x) \Big[4 \log \cosh \frac{x}{2}  +2(\cosh (x)-1) \log \coth \frac{x}{2}  \Big],  \;\;\;\;\; x> 0 \]
  Then $\coth^2 \frac{x}{2} $ has the cross ratio distribution.  The pair $(x_\Sigma,y_\Sigma)$ determines a punctured torus {\em group} $\langle f,g\rangle$ uniquely up to conjugacy.  However Nielson moves on generators give the same quotient but different hyperbolic generators.  For instance the group $\langle fg^3, g \rangle = \langle f, g \rangle$ and both $[f,g]$ and $[fg^3,g]$ are simultaneously parabolic.
 \bigskip
 
 What we would most like to do however is relate the distribution of $[Q]$ to the Teichm\"uller metric on the space of once punctured tori.  To do this we have to describe the relationship between the cross ratio of the vertices of an ideal quadrilateral and it's conformal modulus.  Unfortunately this is not so simple,  in fact Nehari decribes the relationship as ``extremely unobvious" in his book,  \cite[pp 202]{Nehari}.  The particular problem we face of computing the conformal modulus of the ideal quadrilateral $Q$ was also suggested in \cite{KRV}.

\section{Lattices, moduli and the Teichm\"uller metric.} 

We recall that $T^{2}_{*}(\alpha)=  (\IC \setminus \Lambda_{\alpha}(0))/\Lambda_{\alpha}$ as above and have $\alpha>0$.  Rotation and scaling assures us that we may as well assume $\alpha\geq 1$. A fundamental period for the lattice is $[0,\alpha]\times [0,1]$ and this has conformal modulus $K_\alpha= \alpha\geq 1$.  This number is also (by definition) the conformal modulus of $T^{2}_{*}(\alpha)$.  

It is not entirely trivial,  but a bit of reflection should convince one that $K_\alpha$ is the maximal distortion of the extremal quasiconformal mapping between $T^{2}_{*}(1)$ and $T^{2}_{*}(\alpha)$ and further,  the maximal distortion of the extremal quasiconformal mapping between $T^{2}_{*} (\alpha_1)$ and $T^{2}_{*} (\alpha_2)$ is $\max\{|\alpha_1/\alpha_2|,|\alpha_2/\alpha_1|\}$.  Very roughly,  we can define the Teichm\"uller distance between two Riemann surfaces of the same topological type,  say $\Sigma_1$ and $\Sigma_2$,  as 
\[
d_T(\Sigma_1,\Sigma_2) = \inf \{K(f): f:\Sigma_1\to\Sigma_2  \mbox{  is $K$-quasiconformal}\}.
\]
We refer to \cite{A,AIM} for the theory of quasiconformal mappings and for a thorough account of Teichm\" uller theory,  good places to start are \cite{B,G,Hb}.

\medskip

Therefore,  for rectangular punctured tori we would like to relate the conformal modulus of a fundamental period and the cross ratio. A way to do this is to construct a conformal mapping from a fundamental hyperbolic quadrilateral for the group $\langle f,g \rangle$ as above and the period $R_\alpha=[0,\alpha]\times [0,1]$.  

It turns out that this construction is a classical problem investigated separately by Hilbert and Klein, \cite{H,K},  and more recently in \cite{KRV} because of these connections.  We follow the careful exposition of Eremenko,  \cite{Eremenko} (with some notational changes) who considered the problem in relation to  the precise constants in Landau's theorem.

We desire a conformal map $\phi$ of a circular quadrilateral $Q$ (having zero angles, inscribed in the unit circle, symmetric with respect to reflections in coordinate axes) onto a fundamental rectangle $R=R_\alpha$ of the lattice with $\phi(0)$ as the cente of $R$. 

Let ${\cal P}$ be the Weierstrass function of the rectangular lattice with real period 2 and pure imaginary period $2i \tau$. We set
\begin{equation}
{\displaystyle \wp }(z) = \frac{1}{4} ({\cal P}(z +1 +i \tau)-{\cal P}(2)) = \frac{\pi ^2 \vartheta _1^{\prime }\left(0,e^{-\frac{\pi }{\tau}}\right){}^2 \vartheta _1\left(\frac{\pi  z}{2 \tau},e^{-\frac{\pi }{\tau}}\right){}^2}{16 \tau^2 \vartheta _3\left(0,e^{-\frac{\pi }{\tau}}\right){}^2 \vartheta _3\left(\frac{\pi  z}{2 \tau},e^{-\frac{\pi }{\tau}}\right){}^2}, 
\end{equation}
in terms of the Elliptic Theta functions, \cite{Ak}.  Now ${\displaystyle \wp }$ is real on both the real and the imaginary axis, in fact it maps the rectangle $[0,1]\times[0,\tau]$ - one quarter of the fundamental rectangle - onto the lower half-plane.  ${\displaystyle \wp }$ is holomorphic in the closure of  this rectangle except at the point $1+i\tau$ where it has a pole of second order.  Ultimately we will want to solve an equation involving the Schwarzian derivative and this is done by comparing the ratio of two solutions to Lamer's equation.  Thus we consider the eigenvalue problem for complex valued $w$
\begin{equation}\label{lamer}
\frac{\partial^2}{\partial z^2} \; w + {\displaystyle \wp }(z) w = \lambda w.
\end{equation}
Here $\lambda = \lambda(\tau)$ is called an accessory parameter and is to be chosen later.  In fact, identifying $\lambda$ is the main geometric problem, \cite{Nehari2}. Since ${\displaystyle \wp }(z)<0$ for $z\in [0,1]$,  the largest eigenvalue of (\ref{lamer}) with boundary values $w'(0)=w'(1)=0$ is $\lambda_-<0$.  Similarly if we put $z=it$ in (\ref{lamer}) we obtain the equation 
\[ \frac{\partial^2}{\partial t^2} \; w =({\displaystyle \wp }(it)-\lambda) w, \]
whose smallest  eigenvalue with boundary data $w'(0)=w'(i\tau)=0$ is $\lambda_+>0$ since ${\displaystyle \wp }(it)>0$ for $t\in [0,\tau]$.  

Let $c(z),s(z)$ be two solutions of Lamer's equation (\ref{lamer}) with boundary data $c(0)=1,c'(0)=0$ and $s(0)=0,s'(0)=1$.  One can check that $c$ and $s$ are real on $[0,1]$,  while $c$ is real on $[0,i\tau]$ and $s$ is purely imaginary on this interval.  We define the four real quantities $a_j,r_j$,  $j=1,2$ by
\begin{equation}\label{circles}
\left. \begin{array}{lll}
2 a_1 & = &   s(1)/c(1)+s'(1)/c'(1), \\
2 r_1 & = &   s(1)/c(1)-s'(1)/c'(1), \\
2 i a_2 & = &   s(i\tau)/c(i\tau)+s'(i\tau)/c'(i\tau), \\
2 i r_2 & = & s(i\tau)/c(i\tau)+s'(i\tau)/c'(i\tau) .
\end{array} \;\;\;\;\;\; \right\}
\end{equation}
Now if we define
\[ f(z) = \frac{c(z)}{s(z)}, \]
then $f$ solves the Schwarzian derivative equation
\begin{equation}\label{sdeqn}
\frac{f'''}{f'} - \frac{3}{2}\left(\frac{f''}{f'}\right)^2 = 2 ({\displaystyle \wp }-\lambda).
\end{equation}
In fact every solution to (\ref{sdeqn}) is a ratio of linearly independent solutions to (\ref{lamer}).  By construction our $f$ is real on the real axis and purely imaginary on the imaginary axis.  It is also locally univalent.  The Sturm Comparison Theorem shows us that $c(z)$ does not vanish for $\lambda\in [\lambda_-,\lambda_+]$ when $z\in [0,1]\cup [0,i\tau]$ and so $f$ is holomorphic on these intervals.  The reflection principle also  shows $f(\zbar)=\overline{f(z)}$ and $f(-\zbar)=-\overline{f(z)}$.  It now requires a digression to see that the images of $f(\{i\tau\}\times [0,1])$ is a subarc of the circle $\IS(a_2,r_2)$ and $f(\{1\}\times \times [0,i\tau])$ is a subarc of the circle $\IS(a_1,r_1)$ which we leave the reader to find in \cite{Eremenko} or Nehari's book \cite{Nehari}.
  \begin{center}
\scalebox{0.55}{\includegraphics*[viewport=20 520 760 760]{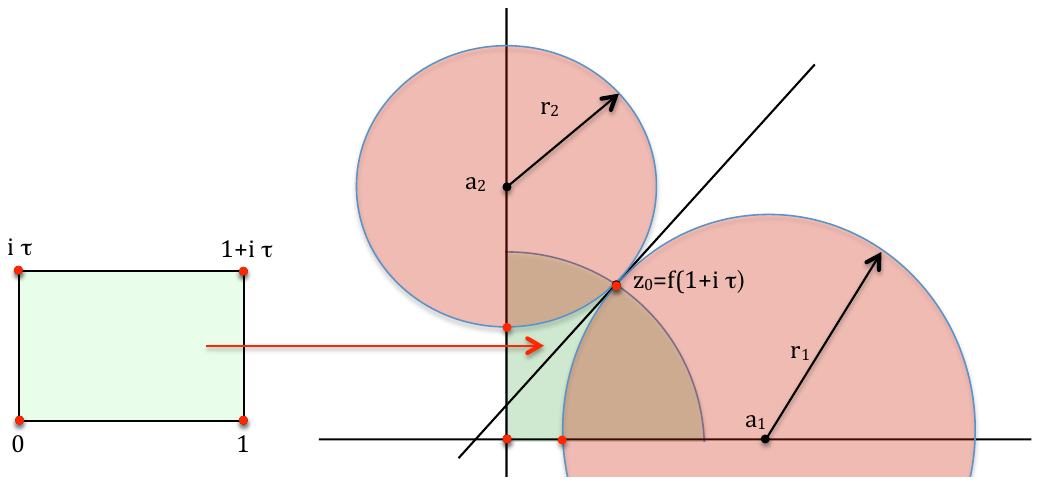}} \\
The mapping $f$.
\end{center}
Set $z_0=f(1+i\tau)$,  the point of tangency of the two circles $\IS(a_1,r_1)$ and $\IS(a_2,r_2)$.  Now if the common tangent line through $z_0$ passes through the origin,  then $\frac{1}{|z_0|} \; f([0,1]\times [0,i\tau])$ is an ideal triangle in the hyperbolic disk $\ID$.  Then developing the obvious reflections across lines of the rectangle and then subarcs of circles shows $f$ to be the local inverse of the universal cover map from $\IC\setminus \Lambda(0)$ to $\ID$.  A little geometric reflection shows that this tangency condition is satisfied when $\sin(\arg(z_0))=r_1/a_1$ and $\cos(\arg(z_0))=r_2/a_2$.  That is when
\begin{equation}\label{test1}
\Big(\frac{r_1}{a_1}\Big)^2+\Big(\frac{r_2}{a_2}\Big)^2 =1.
\end{equation}
Given $\tau$ and the definitions at (\ref{circles}) this becomes a condition on $\lambda$.  In the quadrant that $z_0$ lies, 
as $\lambda\to\lambda_-$, the right circle becomes a vertical line, and when $\lambda\to \lambda_+$ the top circle becomes a horizontal line. The intermediate value theorem implies the existence of the $\lambda$ we seek.  A graph of $\lambda$ - the accessory parameter - is below.
  \begin{center}
\scalebox{0.55}{\includegraphics*[viewport=20 450 660 760]{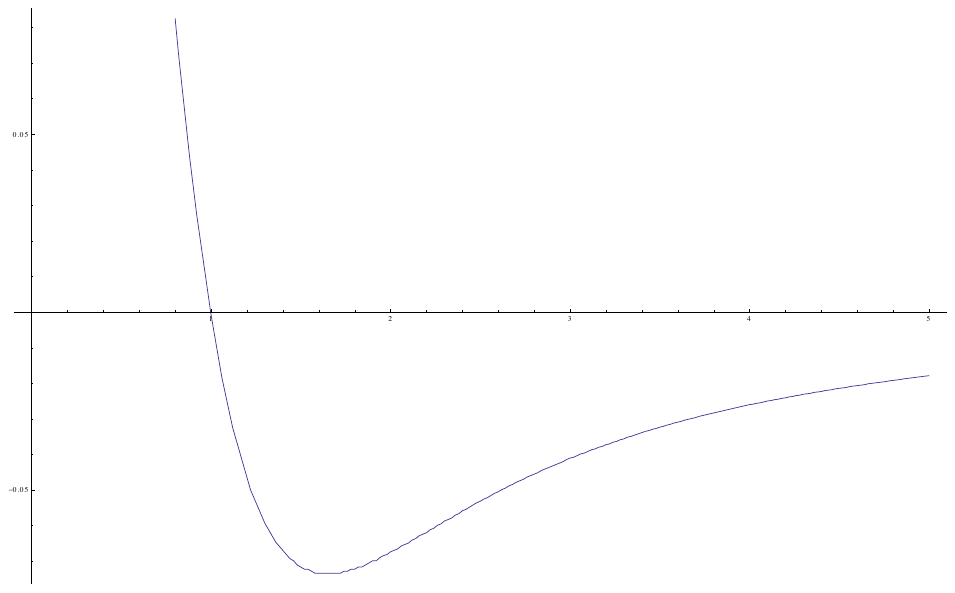}} \\
{\em The parameter $\lambda$ as a function of conformal modulus.}
\end{center}
Actually from a computational point of view the condition (\ref{test1}) is problematic as $1$ is the maximum value of this function.  In practice we use the equivalent condition $\frac{a_{2}^{2}} {\sqrt{a_{1}^{2}+a_{1}^{2}}}-r_1=0$ which crosses the  $y$-axis so root finding algorithms are much faster.  Now when (\ref{test1}) holds a little trigonometry yields
\[
z_0 = \sqrt{a_{1}^{2}-r_{1}^{2}}\Big( \frac{\sqrt{a_{1}^{2}-r_{1}^{2}}}{a_1}+i \;\frac{r_1}{a_1}\Big).
\]
Then the cross ratio
\begin{equation} [z_0,-\overline{z_0},-z_0,\overline{z_0}]= \frac{4|z_0|^2}{ 4 \Re e^2(z_0)} = \frac{a_{1}^{2}}{a_{1}^{2}-r_{1}^{2}} 
= \Big(\frac{a_{2}}{r_{2}}\Big)^2 \geq 1. \end{equation}
We have written code in both Matlab and Mathematica to effect the calculation of these parameters and to draw the graphs and distributions that follow.  The author is happy to distribute these.

\section{Asymptotics.}
In order to discuss approximations and help verify numerical calculations we now discuss the assymptotics.  We first have the following theorem.
\begin{theorem}\label{assympt}  There are constants $c_1,c_2\geq 0$ such that
\begin{equation}
\frac{\pi}{2} \sqrt{[Q]} - c_1 \leq m(Q) \leq \frac{\pi}{2} \sqrt{[Q]} + c_2
\end{equation}
\end{theorem}
Indeed our computational experimentation suggests $c_2=0$ and $c_1=\frac{\pi}{2}$ suffices,  with $c_2\sim 0.9$ for large $[Q]$.

\subsection{Warschawski quadrilaterals.}  We want to establish bounds by comparison of quadrilaterals of a particular form.  The following lemma is an easy consequence of a known modulus estimate of Warschawski.
\begin{lemma}\label{wars}
Let $\alpha \leq 1$ and $\lambda\geq 1$.  Consider the `Warschawski quadrilateral' $Q$ bounded by three line segments $[0,i\alpha/\lambda]$,  $[0,1/\lambda]$, $[1/\lambda, 2i\alpha/\lambda]$ and the graph of $y= \alpha(1+\lambda^2x^2)/\lambda$.  The the modulus of $Q$ is
\begin{equation}
m(Q) = \frac{\pi}{4\alpha}+O(1),  \hskip10pt \mbox{as $\alpha\to 0$}.
\end{equation}
\end{lemma}
\noindent{\bf Proof.} We note that $m(Q)=m(\lambda Q)$ and $\lambda Q$ is the quadrilateral bounded by the three line segments $[0,i\alpha]$, $[0,1]$ and $[1,2i\alpha]$.  We reparameterise the bounding graph to see it has the form $y = \alpha(1+s^2)$,  $s\in[0,1]$. It is  known that the modulus of $\lambda Q$, is equal to
\[m(Q)= \int_0^1\frac{dx}{\alpha(1+x^2)}=\frac{\arctan(1)}{\alpha}=\frac{\pi}{4\alpha}+O(1),  \hskip10pt \mbox{as $\alpha \to 0$}.\]
Here,  the lower estimate is the Ahlfors distortion theorem, the upper
estimate is Warschawski theorem. An elementary proof of
the special case of Warschawski theorem needed for our purposes
can be found in Evgrafov's books \cite{ev1,ev2}.  \hfill $\Box$

\bigskip

\noindent{\bf Proof of Theorem \ref{assympt}.} 
We consider the ideal hyperbolic quadrilateral  $Q_a$  with vertices
$$(a,-\overline{a},-a,\overline{a}),\quad a=\exp(i\theta),\quad\theta\in(0,\pi/4).$$
which we map conformally onto the rectangle $R_a$ with vertices $\{m,m+i,i,0\}$. Now by definition
\begin{equation}
[Q_a]=\frac{1}{\sin^2 \theta} \geq 2,
\end{equation}
and we seek an asymptotic formula for the modulus $m$ of the  rectangle $R_a$ as $\theta\to 0$.

\medskip

As earlier,  it is enough to map  one quarter of our quadrilateral $Q_a$
into the first quadrant since the map extends by symmetry.  Call this quarter in the positive quadrant $Q_{a}^{+}$.  The symmetries of the conformal mapping to $R$ establish that
\[  m(Q_{a}^{+}) = m(R_a/2) = m(R_a). \]
We want to enclose $Q_{a}^{+}$ by two Warschawski quadrilaterals.  

We have the following description of $Q_{a}^{+}$.  The upper side is the arc of a circle orthogonal to the unit circle.
The center of this circle is a point $iy_0$. Let us find it: the tangent
line to the unit circle at $a$ has equation
$t\mapsto e^{i\theta}+ite^{i\theta}$.
This line crosses the imaginary line when $t=\cot\theta$, so
$y_0=\sin\theta+\cos\theta\cot\theta$ and the circle containing the upper side has center $y_0$ and radius
$r=\cot\theta$. We can write the equation of this circle as
\begin{equation} y(x)=\sin \theta  +\cos  \theta   \cot  \theta  -\sqrt{\cot ^2 \theta  -x^2}.
\end{equation}
The other circle, which has a subarc as an edge of $Q_{a}^{+}$, is the circle with radius $\tan \theta$ centered at $1/\cos\theta$.  We now describe two quadrilaterals:
\begin{eqnarray*}
R_0 & \mbox{bounded by} & [0,\frac{1-\sin\theta}{\cos\theta} ],\;\; [0,i\tan(\theta/2)], \;\;[\frac{1-\sin\theta}{\cos\theta},2i\tan(\theta/2) ] 
\\&& \mbox{and the graph $\{y=\tan \left(\frac{\theta }{2}\right) \left(\frac{x^2 \cos ^2(\theta )}{(1-\sin (\theta ))^2}+1\right) \}$}, \\
R_1 & \mbox{bounded by} & [0,\cos\theta],\;\; [0,\frac{i}{2}\sin\theta], \;\;[\cos\theta,i\sin\theta] ,  
\\&& \mbox{and the graph $\{y=\frac{1}{2} \sin (\theta ) \left(\frac{x^2}{\cos ^2(\theta )}+1\right) \}$},  
\end{eqnarray*}
Then Lemma \ref{wars} gives us the two estimates
\[ m(R_0)= \frac{\pi(1-\sin\theta)}{4\tan(\theta/2) \cos(\theta)} +O(1), \;\;\;\;  m(R_1)= \frac{\pi \cos\theta}{2\sin \theta } +O(1) .\]
 \begin{center}
\scalebox{0.65}{\includegraphics*[viewport=30 520 660 800]{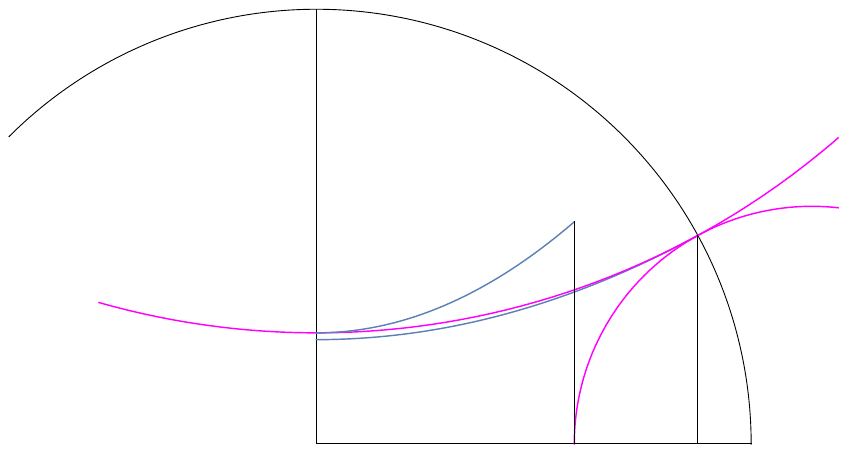}} \\
\end{center}
In order to bound $m(Q_{a}^{+})$ we first need to show that for $\theta$ small and $0\leq x\leq 1$ the graphs of the Warshawki quadrilaterals bound the upper curve of $Q_{a}^{+})$.  That is
\[\frac{1}{2} \sin (\theta ) \Big[ \frac{x^2}{\cos ^2(\theta )}+1\Big] \leq \frac{1}{\sin \theta }-\sqrt{\cot ^2(\theta )-x^2} \leq \tan \frac{\theta }{2} \Big[\frac{x^2 \cos ^2(\theta )}{(1-\sin (\theta ))^2}+1\Big]\]
For our purposes,  given all the functions involved are analytic in larger ranges of $x$ and $\theta$ and $x$ lies in $[0,1]$,  it suffices to observe the power series expansions below and note the common first term.
\begin{eqnarray*}
\frac{1}{2} \sin (\theta ) \left(\frac{x^2}{\cos ^2(\theta )}+1\right)  & \sim & \frac{1}{2} \theta  \left(x^2+1\right)  + \theta ^3 \left(\frac{5 x^2}{12}-\frac{1}{12}\right)+ \cdots \\
\frac{1}{\sin (\theta )}-\sqrt{\cot ^2(\theta )-x^2} & \sim & \frac{1}{2} \theta  \left(x^2+1\right)+\frac{1}{24} \theta ^3 \left(3 x^4+4 x^2+1\right) +\cdots \\
\tan \left(\frac{\theta }{2}\right) \left(\frac{x^2 \cos ^2(\theta )}{(1-\sin (\theta ))^2}+1\right) & \sim &\frac{1}{2} \theta  \left(x^2+1\right)+\theta ^2 x^2+ \frac{1}{24}\theta ^3 \left(25 x^2 +1\right)+\cdots
\end{eqnarray*} 
One could of course show the differences between these functions either increases with $x$, or decreases with $x$,  using calculus,  but another virtue of the series above is that they quickly reveal
\[ \frac{\pi(1-\sin\theta)}{4\tan(\theta/2) \cos(\theta)} = \frac{\pi}{2\sin\theta} -\frac{\pi}{2} + O(\theta) ,   \]
and 
\[ \frac{\pi \cos\theta}{2\sin \theta } =  \frac{\pi}{2\sin\theta} +O(\theta), \]
and hence
\[ m(R_0)= \frac{\pi \sqrt{[Q_a]}}{2} +O(1), \;\;\;\;\; {\rm and}\;\;\;\;\; m(R_1)= \frac{\pi \sqrt{[Q_a]}}{2} +O(1). \]
All that remains to complete the proof is to make some elementary observations comparing the moduli of the curve families associated to the three quadrilaterals $R_0,R_1$, and $Q_{a}^{+}$.  For example,  if $\rho(z)$ is an admissible function for the family of curves joining the left and right side of $R_0$ (within $R_0$),  then the $\rho$-length of any curve joining the left and right sides of $Q_{a}^{+}$ (withing $Q_{a}^{+}$) is at least one.  Thus $\chi_{Q_{a}^{+}}\rho$ is admissible to $Q_{a}^{+}$.  It follows $m(Q_{a}^{+})\geq m(R_0)$.  A similar argument is made for $R_1$ and the result follows. \hfill $\Box$.
  
 \begin{center}
\scalebox{0.65}{\includegraphics*[viewport=30 450 760 760]{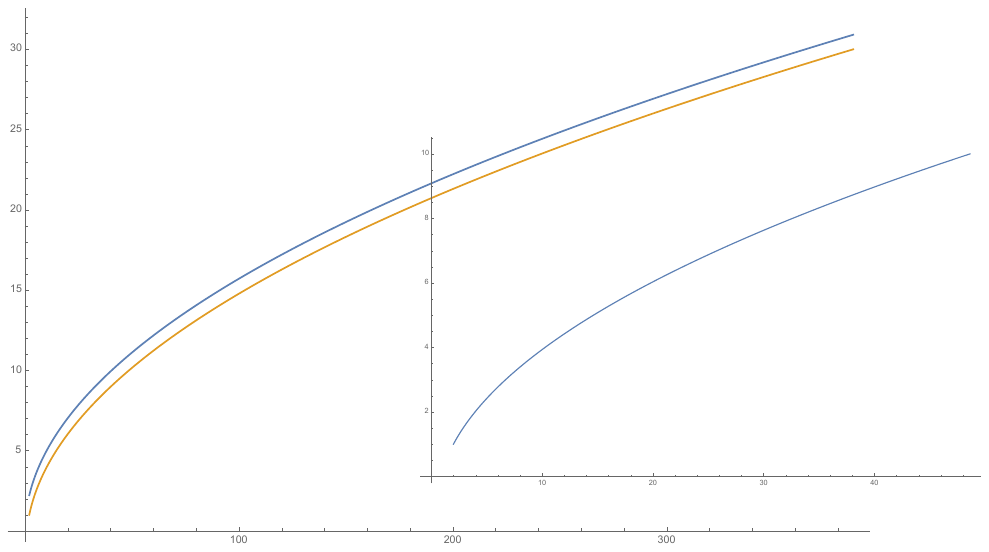}} \\
{\em The graph of ${\cal CR}$. $x$-axis: cross ratio,  $y$-axis: modulus.  Lower curve is the approximation $m(Q) = \frac{\pi}{2}\sqrt{[Q]}$.  Inset detail.}
\end{center}

\section{Distributions.}
We have now defined an increasing homeomorphism ${\cal CR}:[1,\infty]\to[2,\infty]$ which takes a modulus $m\geq 1$ to a cross ratio $[Q]\geq 2$.  Examining the symmetries of modulus and of the cross ratio we can extend ${\cal CR}$ to $[0,\infty)\to[0,\infty)$ as an analytic function by the rule
\begin{equation}\label{crfeqn}
{\cal CR}(1/m) = \frac{{\cal CR}(m)}{{\cal CR}(m)-1}, \hskip10pt {\cal CR}(1)=2.
\end{equation}
Restricting to random cross ratios bigger than $2$ gives us the p.d.f.  
\begin{equation}
{\cal F}(r) = \frac{6}{\pi^2}\left( \frac{  \log (r)}{(r-1) r} + \frac{ \log \left( \frac{r}{r-1}\right)}{r}  \right), \hskip10pt \mbox{if $r\geq 2$}.
\end{equation}
Notice that ${\cal F}$ satisfies the functional equation
\[ {\cal F}(r)(r^2-1) = {\cal F}(r/(r-1)),\hskip10pt r>1,  \]
and thus ${\cal F}({\cal CR}(1/m)) = {\cal F}({\cal CR}(m) )({\cal CR}^2(m) -1)$.

We can now compute the probability distribution function for the modulus as 
\begin{equation}
{\cal M}(m) = {\cal F}({\cal CR}(m)) {\cal CR}'(m) ,\hskip10pt m\geq 1  .
\end{equation}
Unfortunately we cannot directly compute ${\cal CR}'(m)$.  Differentiating (\ref{crfeqn}) yields
\begin{eqnarray*}
{\cal CR}'(1/m) & = &  \frac{{m^2\cal CR}'(m)}{({\cal CR}(m)-1)^2} ,\\
{\cal CR}''(1/m) & =& -\frac{m^4 {\cal CR}''(m)}{({\cal CR}(m)-1)^2}+\frac{2 m^4{\cal CR}'(m)^2}{({\cal CR}(m)-1)^3}-\frac{2 m^3 {\cal CR}'(m)}{({\cal CR}(m)-1)^2}, \\
&\vdots & .
\end{eqnarray*}
The second equation here shows us that 
\begin{eqnarray*} 
{\cal CR}''(1) & =&   {\cal CR}'(1)^2 -    {\cal CR}'(1) ,
\end{eqnarray*}
which gives us the series near $m=1$,
\begin{equation}\label{taylor} {\cal CR}(m) = 2+a(m-1)+(a^2-a)(m-1)^2+ \cdots,  \hskip10pt a= {\cal CR}'(1). \end{equation}
Curiously the functional equation (\ref{crfeqn})  does not  determine the higher odd degree coefficients of the Taylor series expansion at $m=1$.
 \begin{center}
\scalebox{0.7}{\includegraphics*[viewport=30 470 760 760]{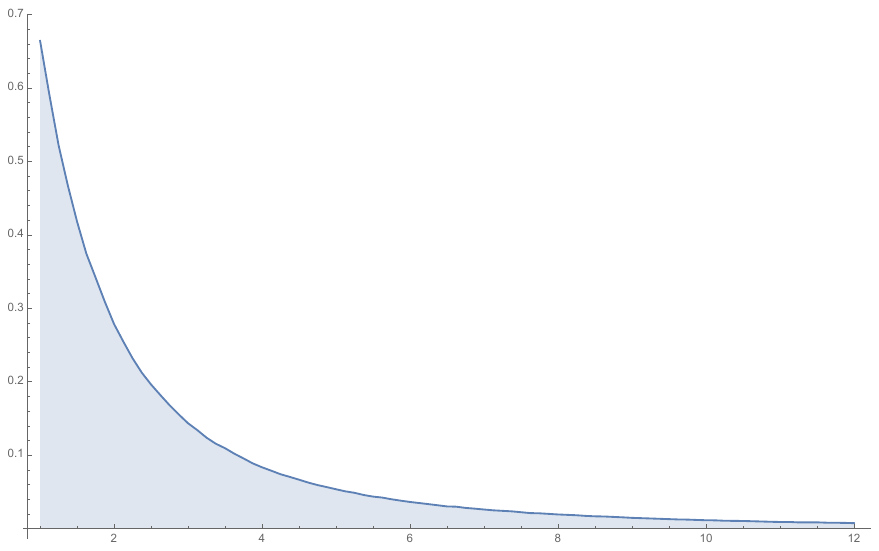}} \\
{\em The p.d.f. ${\cal M}(m)$ for the modulus of a rectangular punctured torus.  }
\end{center}
For large $m$ we have ${\cal M}(m)\sim \frac{192 \log m}{\pi ^5 m^3}$.  Using (\ref{taylor}) we can compute the first few terms of the series for ${\cal M}(m)$ for $ m\sim 1$.
\begin{eqnarray*}\lefteqn{{\cal M}(m)}\\
& = & \frac{a \log (64)}{\pi ^2}-\frac{6 (2-a) a \log (2)}{\pi ^2}(m-1) -\frac{3  a^2 (a+8 (2 a-3) \log (2)) }{4 \pi ^2} (m-1)^2 \\&& +\frac{3 a^2  (a (a (\log (16)-1)+2)-8 \log (2))}{2 \pi ^2}(m-1)^3+O\left((m-1)^4\right).
\end{eqnarray*}
\medskip
Our computation experiments suggest the conjectural value $a=\frac{\pi}{2}$.  

\bigskip

The logarithmic transformation of this p.d.f. will give us the p.d.f. ${\cal T}(d)$ for the Teichm\"uller distance of a rectangular punctured torus $T^{2}_{*}(m)$ of modulus $m$ to the square punctured torus $T^{1}_{*}(1)$.  The change of variables formula for probability distributions gives us the formula
\[ {\cal T}(d) = {\cal M}(e^d)e^d, \hskip15pt d>0. \]
Notice that as soon as ${\cal M}'(1)+{\cal M}(1)  >0$ this distribution is initially increasing.  This condition is equivalent to $a>1$,  a condition which is verifiable computationally.
 \begin{center}
\scalebox{0.7}{\includegraphics*[viewport=30 500 760 760]{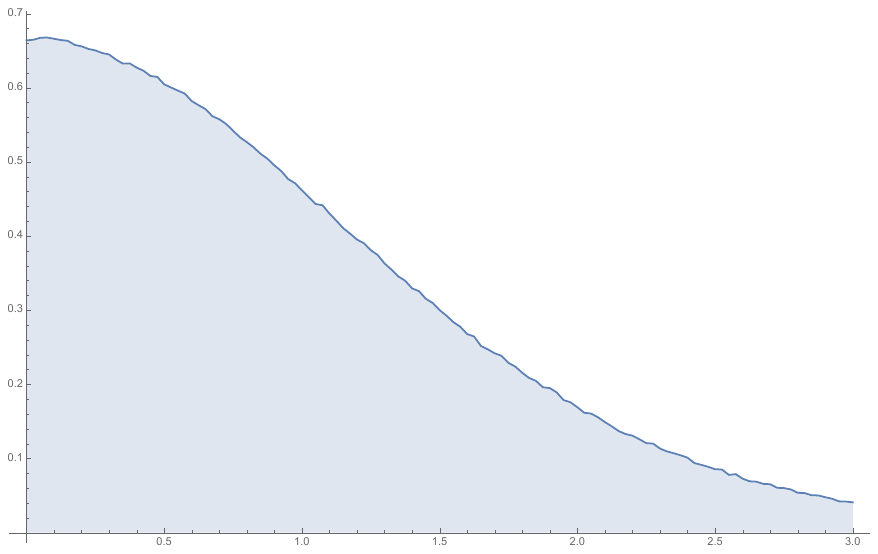}} \\
{\em The p.d.f. for the random variable ${\cal T}(d)$ : Teichm\"uller distance of a random punctured torus to the square punctured torus.}
\end{center}
For $d$ small we have the expansion (with the conjectural value $a={\cal CR}'(1)$), 
\begin{eqnarray*}
{\cal T}(d) & = & \frac{a \log (64)}{\pi ^2}+\frac{6 (a-1) a   \log 2}{\pi ^2} d -\frac{3  a \left(a^2+4 (a-1) (4 a-5) \log (2)\right)}{4 \pi ^2}d^2 \\ && +\frac{(a-1) a  \left((6 (a-3) a+13) \log (4)-3 a^2\right)}{2 \pi ^2} d^3+O\left(d^4\right) .
\end{eqnarray*}
When we put in $a=\frac{\pi}{2}$ we have
\begin{eqnarray*}
{\cal T}(d) & = & \frac{3 \log 2}{\pi}\Big[1+  \left(\frac{\pi}{2}-1\right) d\Big] - O(d^2).
\end{eqnarray*}
Using the approximation $m=\frac{\pi}{2}\sqrt{[Q]}$ we find that ${\cal T}(d)\sim \frac{192  d}{\pi ^5 e^{2 d}}$ for large $d$.  Thus this distribution has all moments. 
Using various approximations to the function ${\cal CR}(m)$ we find expected value is roughly $1$ with median $\sim 0.779$ and standard deviation $0.803$.

\section{Quasi-M\"obius mappings.} We close with a few minor observations on the distortion of cross ratios under the induced boundary map for quasiconformal mappings of the disk.   A homeomorphism $f_0:\IS\to\IS$ is said to be quasim\"obius if there is an increasing function $\eta:[0,\infty)\to[0,\infty)$ such that for all quadruples $\{z_0,z_1,z_2,z_3\}$ we have
\begin{equation}
f_0[z_0,z_1,z_2,z_3] = [f_0(z_0),f_0(z_1),f_0(z_2),f_0(z_3)] \leq \eta([z_0,z_1,z_2,z_3]).
\end{equation}
The function is weakly $M$-quasim\"obius if whenever $[z_0,z_1,z_2,z_3]=2$ we have $f_0[z_0,z_1,z_2,z_3]\leq 2 M$.
\begin{theorem}  Let $z_0,z_1,z_2,z_3\in \IS$ and $z_0',z_1',z_2',z_3'\in \IS$ be two sets of anticlockwise ordered points in the circle.  Then there is a $K$-quasiconformal mapping $f:\ID\to\ID$ with quasim\"obius boundary values $f_0:\IS\to\IS$ such that $f_0(z_i)=z_i'$,  $i=0,\ldots,3$ and
\begin{equation}
K =\max \left\{ \frac{ {\cal CR}^{-1}([z_0,z_1,z_2,z_3]) }{{\cal CR}^{-1}([z_0',z_1',z_2',z_3'])}, \;\;  \frac{{\cal CR}^{-1}([z_0',z_1',z_2',z_3'])}{{\cal CR}^{-1}([z_0,z_1,z_2,z_3]) }\right\}.
\end{equation}
\end{theorem}
If we fix $[z_0,z_1,z_2,z_3]=2$,  then our asymptotics identified above yield the following information on the weak quasim\"obius constant:
\[  M \sim \frac{\pi}{4} \sqrt{[z_0',z_1',z_2',z_3']}, \;\;\;\mbox{when $[z_0',z_1',z_2',z_3'] >> 2$} ,\]
and 
\[  M \sim  1+\frac{2 }{\pi} \;\epsilon, \;\;\;\mbox{when $2\leq [z_0',z_1',z_2',z_3'] \leq 2+\epsilon$,  $\epsilon <<1$.} \]
There is another obvious map here.  Namely if $[z_0,z_1,z_2,z_3]=2$ then we map to the upper half-space so the vertices are $-1,0,1,\infty$.  Then the $K$ quasiconformal mapping 
\begin{equation}
x+i y \in \IH^2 \mapsto \left\{ \begin{array}{ll} x+i y & \mbox{ if $x\leq 0$ } \\  Kx+i y & \mbox{ if $x\geq 0$ }\end{array}  \right.,
\end{equation}
gives image cross ratio $[-1,0,K,\infty] = K+1$,  so the above estimates are far better.

\medskip

eremenko@math.purdue.edu\\

g.j.martin@massey.ac.nz

\end{document}